\documentclass[11pt,oneside,leqno]{amsart}
\usepackage{amsxtra}
\usepackage{amsopn}
\usepackage{amsmath,amsthm,amssymb}
\usepackage{amscd}
\usepackage{color}
\usepackage{amsfonts}
\usepackage{latexsym}
\usepackage{verbatim}
\usepackage{longtable}
\usepackage{pdflscape}
\usepackage{graphicx}
\usepackage{pb-diagram}

\theoremstyle{plain}
\newtheorem{theorem}{Theorem}[section]
\newtheorem*{theorem*}{Theorem}

\newtheorem{prop}[theorem]{Proposition}
\newtheorem{cor}[theorem]{Corollary}
\newtheorem{rem}[theorem]{Remark}

\newtheorem*{mt*}{Main Theorem}
\newcommand\g{{\mathfrak{g}}}
\newcommand\h{{\mathfrak{h}}}

\newcommand\m{{\mathfrak{m}}}

\newcommand\R{{\mathbb R}}

\newcommand\W{\mathcal W}

\begin{document}
\title[]{Four-dimensional pseudo-Riemannian \\ homogeneous Ricci solitons}
\author{Giovanni Calvaruso and Anna Fino}
\date{}

\address{Giovanni Calvaruso: Dipartimento di Matematica e Fisica\lq\lq E. De Giorgi\rq\rq \\
Universit\`a del Salento\\
Prov. Lecce-Arnesano \\
73100 Lecce\\ Italy} 
\email{giovanni.calvaruso@unisalento.it}

\address{Anna Fino: Dipartimento di Matematica \\ Universit\`a di Torino\\
Via Carlo Alberto 10 \\
10123 Torino\\ Italy} \email{annamaria.fino@unito.it}

\thanks{This work was partially supported by  MIUR, GNSAGA (Italy)}

\begin{abstract} 
We consider four-dimensional homogeneous pseudo-Riemannian manifolds with non-trivial isotropy  and completely classify the cases giving rise to non-trivial homogeneous Ricci solitons. In particular, we show the existence of non-compact homogeneous (and also invariant) pseudo-Riemannian 
Ricci solitons which are not isometric to solvmanifolds, and of conformally flat homogeneous pseudo-Riemannian Ricci solitons which are not symmetric.
\end{abstract} 

\maketitle

\section{Introduction}

Ricci solitons were introduced by Hamilton  \cite{Ham} and they are a natural generalization of   Einstein metrics. A pseudo-Riemannian metric  $g$  on a smooth manifold  $M$  is called a  \emph{Ricci soliton}     if there exists a a smooth vector field $X$, such that
\begin{equation}\label{solit}
\mathcal{L}_Xg+\varrho=\lambda g,
\end{equation}
where $\mathcal{L}_X$ denotes the Lie derivative in the direction of
$X,$ $\varrho$ denotes the Ricci tensor and $\lambda$ is a real number.
A Ricci soliton $g$  is said to be a \emph{shrinking}, \emph{steady} or
\emph{expanding} according to whether $\lambda>0$, $\lambda=0$ or $\lambda<0$,
respectively.

Ricci solitons are  the self-similar solutions of the  
\emph{Ricci flow} and are important in understanding its singularities. A survey and further references on the geometry of Ricci solitons may be found in \cite{Cao}. The interest  in Ricci solitons has also risen among theoretical physicists in
relation with String Theory \cite{AW},\cite{Fr},\cite{Kh}. After their introduction in the Riemannian case, the study of pseudo-Riemannian Ricci solitons attracted a growing number of authors (see for instance \cite{isr}-\cite{CF},\cite{Cas},\cite{PT}). 

{If $M=G/H$ is a homogeneous space, a {\em homogeneous Ricci soliton} on $M$ is a $G$-invariant metric $g$ for which equation \eqref{solit} holds. In particular, by an {\em invariant Ricci soliton} we mean a homogeneous one, such that equation \eqref{solit} is satisfied by an invariant vector field. It is a natural question to determine which homogeneous manifolds $G/H$ admit a $G$-invariant Ricci soliton \cite{L}. }  All known examples of homogeneous Riemannian Ricci soliton metrics on non-compact homogeneous manifolds are isometric to some solvsolitons, that is, to invariant Ricci solitons on a solvable Lie group  (\cite[Remark 1.5]{Jab}).

The difference between Riemannian and pseudo-Riemannian settings lead to different results concerning the existence of homogeneous Ricci solitons. In fact, {although there exist three-dimensional Riemannian homogeneous Ricci solitons \cite{BD},\cite{L},} there are no left-invariant Ricci solitons on three-dimensional Riemannian Lie groups  \cite{Cerbo},  while left-invariant Ricci solitons on three-dimensional Lorentzian Lie groups were classified in \cite{isr}. Four-dimensional Ricci solitons on non-reductive homogeneous pseudo-Riemannian manifolds (which do not have a Riemannian countepart, as a homogeneous Riemannian manifold is necessarily reductive) were classified in \cite{CF}.

In the present article we provide a full classification of four-dimensional homogeneous pseudo-Riemannian Ricci solitons in all the cases with non-trivial isotropy. The starting point is the complete local classification of four-dimensional homogeneous pseudo-Riemannian manifolds with non-trivial isotropy obtained in \cite{K0} (see also \cite{K}), leading to the remarkable number of $186$ different forms of these spaces. Each of them admits a family of invariant pseudo-Riemannian metrics, depending of a number of real parameters varying from $1$ to $4$. Among them, we were able to determine { $44$} different examples of homogeneous spaces $M=G/H$ for which equation \eqref{solit} holds for some vector fields $X \in \m$ and some invariant metrics which are not Einstein. In particular, for { $41$} of these examples, equation \eqref{solit} is satisfied by an invariant vector field. Some of these examples show that there exist homogeneous (and also invariant) pseudo-Riemannian Ricci solitons which are not isometric to solvsolitons.

In the Riemannian case, most of the known examples of non-trivial Ricci solitons  are  K\"ahler-Ricci solitons,  with  the exception of the homogeneous solitons on nilpotent Lie groups \cite{La2}, the rotationally symmetric Bryant solitons  on $\R^n$ ($n > 2$), Ivey's generalization of these solutions \cite{Iv} and the complete steady gradient Ricci solitons constructed in \cite{DW}. By using the results  in \cite{CF2} we show that {none} of the four-dimensional  non-trivial pseudo-Riemannian homogeneous Ricci solitons is K\"ahler.

Conformally flat Einstein pseudo-Riemannian manifolds have constant sectional curvature. In particular, they are symmetric. Conformally flat homogeneous Riemannian manifolds are always symmetric \cite{Ta}. On the other hand, some of our examples show the existence of conformally flat homogeneous pseudo-Riemannian Ricci solitons which are not symmetric.

The paper is organized in the following way. In Section~2, we shall report some basic facts on four-dimensional homogeneous pseudo-Riemannian manifolds with non-trivial isotropy. The classification of homogeneous solutions to \eqref{solit} will be then presented in Section~3.  In Section~4 we shall investigate several geometric properties of four-dimensional pseudo-Riemannian homogeneous Ricci solitons (with non-trivial isotropy). In particular, we show the {non-existence} of  pseudo-K\"ahler examples, i.e., solutions corresponding to holomorphic vector fields $X \in \m$, and we classify the conformally flat examples.   

\medskip
{\noindent \textbf{Acknowledgements.} We would like to thank  Maria Buzano, Eduardo Garcia-Rio, Amirhesam Zaeim and the anonymous referee for  useful  comments on the paper. }
\smallskip

\bigskip \noindent
 \section{Preliminaries}
 \setcounter{equation}{0}

Let $M$ be a homogeneous space, $H = H_x$ the stabilizer of an arbitrary
point $x$ of $M$, and $(\g, \h)$ the pair of Lie algebras corresponding to the
pair $(G,H)$ of Lie groups. As remarked in \cite{K}, the pair $(\g ,\h)$ locally uniquely defines
the homogeneous space.

With regard to four-dimensional pseudo-Riemannian homogeneous spaces, the classification obtained by Komrakov \cite{K0,K} starts with the list of subalgebras in the Lie algebras $\mathfrak{so}(4)$, $\mathfrak{so}(3, 1)$ and $\mathfrak{so}(2, 2)$.
Up to conjugacy, Lie algebras $\mathfrak{so}(4)$, $\mathfrak{so}(3, 1)$ and $\mathfrak{so}(2, 2)$ respectively admit $6$, $14$ and $28$ distinct Lie subalgebras (some of them depending on parameters). 

Given a subalgebra $\h$ of either $\mathfrak{so}(4)$, $\mathfrak{so}(3, 1)$ or $\mathfrak{so}(2, 2)$, assuming that
$\h$ acts naturally on $\mathbb R^4$, there is a one-to-one correspondence between  faithful generalized
modules $(\h , \mathbb R^4)$ and  corresponding subalgebras of $\mathfrak{so}(p, q)$. This leads to determine $186$ distinct models of four-dimensional homogeneous pseudo-Riemannian manifolds with non-trivial isotropy. Following the notation and the classification used in \cite{K}, the space identified by the type ${\bf n.m^k.q}$ is the one corresponding to the ${\bf q}$-th pair  $(\g, \h)$ of type ${\bf n.m^k}$, where ${\bf n}=\dim (\h)$ ($=1,..,6$), ${\bf m}$ is the number of the complex subalgebra $\h ^{\mathbb C}$ of $\mathfrak{so}(4,\mathbb C)$ and ${\bf k}$ is the number of the real form of $\h ^{\mathbb C}$. When the index ${\bf q}$ is omitted, one refers simultaneously to all homogeneous spaces corresponding to pairs  $(\g, \h)$ of type ${\bf n.m^k}$.

Consider now an arbitrary pair $(\g,\h)$ (that is, the corresponding homogeneous space $M=G/H$, with $H$ connected). Let $\m = \g /\h$ denote the four-dimensional factor space, which identifies with a subspace of $\g$ complementary to $\h$. The different cases are explicitly described in \cite{K}, listing a basis $\{e_j\}$ of $\h$, a basis $\{u_i\}$ of $\m$ and their Lie brackets. The pair $(\g, \h)$ uniquely defines the isotropy representation
$${\psi}: \h \to \mathfrak{gl}(\m), \qquad \psi (x)(y) = [x,y]_{\m} \quad \text{for all} \quad x \in \h, y \in \m.$$ 
A bilinear form $B$ on $\m$ is invariant if and only if ${\psi} (x)^t \circ B +B \circ {\psi}(x)=0$, for all $x \in \h$, where ${\psi} (x)^t$ denotes the transpose of ${\psi} (x)$. In particular, requiring that $B=g$ is symmetric and nondegenerate, this leads to the classification of all invariant pseudo-Riemannian metrics on $G/H$.

With regard to curvature properties, following \cite{K}, an invariant nondegenerate symmetric bilinear form $g$ on $\m$ uniquely defines its invariant linear Levi-Civita connection  $\nabla$, described in terms of the corresponding homomorphism of $\h$-modules $\Lambda : \g \to \mathfrak{gl}(\m)$, such that $\Lambda (x)(y_{\m}) = [x, y]_{\m}$ for all $x\in \h, y \in \g$. Explicitly, one has
\begin{equation}\label{LC}
\begin{array}{l}
\Lambda (x)(y_{\m}) = \frac{1}{2}[x,y]_{\m} +v(x,y), \qquad \text{for all} \; x,y \in \g ,
\end{array}
\end{equation}
where $v: \g \times \g \to \m $ is the $\h$-invariant symmetric mapping uniquely determined by 
$$2g (v(x, y), z_{\m}) = g(x_{m}, [z, y]_{\m}) + g(y_{\m}, [z, x]_{\m}), \qquad \text{for all} \; x, y, z \in \g.
$$
The curvature tensor is then determined by the mapping
$R:\m \times \m \to \mathfrak{gl}(\m)$, such that 
\begin{equation}\label{Curv}
R(x, y) = [\Lambda(x), \Lambda(y)]-\Lambda([x, y]),
\end{equation} 
for all $x, y \in \m$. Finally, the Ricci tensor $\varrho$ of $g$, described in terms of its components with respect to $\{u_i\}$, is given by 
\begin{equation}\label{ric}
\varrho (u_i, u_j) = \sum _{r=1} ^4 R_{ri} (u_r,u_j), \qquad i,j=1,..,4.
\end{equation}
It must be noted that whenever $X=\sum x_k u_k \in \m$, equation \eqref{solit} reads as a system of algebraic equations for the components $x_k$ of $X$, namely,
\begin{equation}\label{algsol}
\sum _{k=1} ^4 x_k \big( g([u_k,u_i],u_j)+ g(u_i,[u_k,u_j])\big) +\varrho (u_i,u_j) = \lambda g_{ij}, \quad i,j=1,..,4. 
\end{equation}

\smallskip

{
We shall now describe {only the} homogeneous pseudo-Riemannian four-manifolds $G/H$ which will appear in the classification of non-trivial Ricci solitons. For each of the different cases, numbered as in \cite{K}, we shall report the nonvanishing Lie brackets with respect to the  basis $\{u_j \}$ of $\m$ and the  basis $\{e_i \}$ of $\h$, the generic  invariant pseudo-Riemannian metrics $g$ and the corresponding Ricci tensor $\varrho$ with respect to $\{u_j\}$.  Denoted by $\{\theta_j\}$ the basis of $1$-forms dual to $\{u_j\}$, we shall use the following notations: $\theta_i\circ\theta_j= \frac{1}{2}\left(\theta_i\otimes\theta_j+\theta_j\otimes\theta_i\right)$, $\theta_i ^2= \theta_i\circ\theta_i$, for all indices $i,j$. The Ricci tensor with respect to the basis $\{u_i\}$ has been computed by using \eqref{ric} and {\em Maple 16$^\copyright$}.  }

\medskip
{
{\small
${\bf 1.1^1)}$  $\quad  [e_1,u_1]= u_1, \;  [e_1,u_3]= -u_3$,  so that the invariant metrics are given by
$$g=a\theta_1\circ \theta_3+b\theta_2^2+c\theta_2 \circ \theta_4+d\theta_4^2, \quad a(bd-c^2) \neq 0.$$
The possible cases are the following.
$$\begin{array}{cl}
\vphantom{\displaystyle{A^{A^A}}} {\bf 1.1^1.1}: &    [u_1,u_3]=u_2, \;  [u_2,u_4]= u_2,  \; [u_3,u_4]= u_3;  \\[4pt]
&  \varrho=-\frac{b(2a^2+c^2-bd)}{2a(bd-c^2)} \theta_1\circ \theta_3-\frac{b^2 (4a^2 + bd - c^2)}{2 a^2 (bd-c^2)}\theta_2^2-\frac{b c (4a^2 + bd - c^2)}{2 a^2 (bd-c^2)}\theta_2 \circ \theta_4-\frac{c^2 a^2 + bd c^2 - c^4 + 3 b a^2  d}{2 a^2 (bd-c^2)}\theta_4^2.
\vspace{3pt}
\end{array}$$ 
$\begin{array}{cl}
\vphantom{\displaystyle{A^{A^A}}} {\bf 1.1^1.2}: & [u_2,u_4]= p u_2, \; [u_3,u_4]= u_3, \, p \in \R ; \\[4pt]  &  \varrho=-\frac{ab(p+1)}{2(bd-c^2)} \theta_1\circ \theta_3--\frac{b^2 p(p+1)}{bd-c^2}\theta_2^2-\frac{bc p(p+1)}{bd-c^2}\theta_2 \circ \theta_4-\frac{2c^2 p+bd-c^2+2bd p^2}{2(bd-c^2)} \theta_4^2
\vspace{3pt}
\end{array}$ 
$\begin{array}{cl}
\vphantom{\displaystyle{A^{A^A}}} {\bf 1.1^1.5}: &  [u_1,u_3]=  e_1,\; [u_2,u_4]=  u_2 ;  \\[4pt]
&  \varrho=- \theta_1\circ \theta_3-\frac{b^2}{bd-c^2}\theta_2^2-\frac{bc}{bd-c^2}\theta_2 \circ \theta_4-\frac{bd}{bd-c^2} \theta_4^2
\vspace{3pt}
\end{array}$ 
}}

\medskip
{
{\small
${\bf 1.1^2)}$  $\quad  [e_1,u_1]= u_3, [e_1,u_3]= -u_1$,  so that the invariant metrics are given by
\begin{equation}\label{first}
g=a\theta_1^2+b\theta_2^2+c\theta_2 \circ \theta_4+a\theta_3^2+d\theta_4^2, \quad a(bd-c^2) \neq 0.
\end{equation}
The possible cases are the following. \\
$\begin{array}{cl}
\vphantom{\displaystyle{A^{A^A}}} {\bf 1.1^2.1}: &   [u_1,u_3]= -u_2, [u_1,u_4]= u_1, [u_2,u_4]= 2 u_2, [u_3,u_4]= u_3 ;  \\[4pt]
&  \varrho=-\frac{(8a^2+bc-d^2)b}{2a(bc-d^2)} \theta_1^2 -\frac{(16a^2-bc+d^2)b^2}{2a^2(bc-d^2)}\theta_2^2-\frac{(16a^2-bc+d^2)bd}{2a^2(bc-d^2)}\theta_2 \circ \theta_4-\frac{(8a^2+bc-d^2)b}{2a(bc-d^2)}\theta_3^2\\[4pt] & \quad -\frac{12a^2bc+4a^2d^2 -bcd^2+d^4}{2a^2(bc-d^2)}\theta_4^2.
\vspace{3pt}\\[4pt]
\vphantom{\displaystyle{A^{A^A}}} {\bf 1.1^2.2}: &  [u_1,u_4]= u_1, \; [u_2,u_4]= p u_2, \; [u_3,u_4]= u_3, \; p\in \R  ;  \\[4pt]%
& \varrho=-\frac{ab(p+2)}{bc-d^2} \theta_1^2 -\frac{b^2 p(p+2)}{bc-d^2}\theta_2^2-\frac{bdp(p+2)}{bc-d^2}\theta_2 \circ \theta_4-\frac{ab(p+2)}{2bc-d^2}\theta_3^2 -\frac{bc(p^2+2)+2d^2(p-1) }{bc-d^2}\theta_4^2.
\end{array}$ 
$\begin{array}{l}
\\
\vphantom{\displaystyle{A^{A^A}}} {\bf 1.1^2.6},{\bf 1.1^2.7}: \quad [u_1,u_3]= \varepsilon e_1, \; [u_2,u_4]= u_2;  \\[4pt]
\hspace{27mm} \varrho=\varepsilon\theta_1^2 -\frac{b^2 }{(bc-d^2)^2}\theta_2^2-\frac{bd}{(bc-d^2)}\theta_2 \circ \theta_4+\varepsilon\theta_3^2 -\frac{bc}{(bc-d^2)}\theta_4^2, 
\vspace{3pt}
\end{array}$ 
\\[3pt]
where $\varepsilon=1$ for ${\bf 1.1^2.6}$ and $\varepsilon=-1$ for ${\bf 1.1^2.7}$.

\medskip
{
{\small
${\bf 1.3^1)}$  $\quad [e_1,u_3]= u_1, \; [e_1,u_4]= u_2$,  so that the invariant metrics are given by
$$g=a\theta_1\circ \theta_4-a\theta_2 \circ \theta_3+b\theta_3^2+c\theta_3 \circ \theta_4+d\theta_4^2, \quad a \neq 0.$$
The possible cases are the following. \\
$\begin{array}{cl}
\vphantom{\displaystyle{A^{A^A}}} {\bf 1.3^1.1}: &  [u_1,u_2]= - \frac{1}{2} u_2, \; [u_1,u_3]=  u_1 \; [u_1,u_4]=  \frac 1 2 u_4, \; [u_2,u_3]=  \frac 12 u_4; \\[4pt]
&  \varrho=\frac{3d}{8a} \theta_1\circ \theta_4-\frac{3d}{8a} \theta_2\circ \theta_3 +\frac{3(6bd-5c^2)}{8a^2}\theta_3^2   +\frac{3cd}{8a^2} \theta_3 \circ \theta_4+\frac{3d^2}{8a^2}\theta_4^2.
\vspace{3pt}
\end{array}$ 
$\begin{array}{cl}
\vphantom{\displaystyle{A^{A^A}}} {\bf 1.3^1.2}: &  [u_1,u_3]=-p e_1+ (p+1) u_1 +pu_2,  [u_2,u_4]=u_2, \; |p|\leq 1; \\[4pt]
&  \varrho=-\frac 12(p^2+1)\theta_3^2 +\frac 12(p+1) \theta_3 \circ \theta_4-\frac 12\theta_4^2.
\vspace{3pt}
\end{array}$\\ 
$\begin{array}{cl}
\vphantom{\displaystyle{A^{A^A}}} {\bf 1.3^1.3}: &   [u_1,u_3]= u_1 , \; [u_2,u_4]=u_2, \;  [u_3,u_4]= e_1; 
\\[4pt]
&  \varrho=-\frac 12\theta_3^2 +\frac 12 \theta_3 \circ \theta_4-\frac 12\theta_4^2
\vspace{3pt}
\end{array}$ 
\\$\begin{array}{cl}
\vphantom{\displaystyle{A^{A^A}}} {\bf 1.3^1.5}: & [u_1,u_3]= -\frac{p^2 +q}{q -1}e_1+\frac{1+p^2}{q -1} u_2 , \; [u_1,u_4]=  pe_1 + u_1+p u_2, \; [u_2,u_3]=  pe_1 + u_1+p u_2, \\[4pt] &   
 [u_2,u_4]=-q  e_1+(q +1)  u_2, \quad  p \geq 0, \; q \neq 1; \\[4pt]
&  \varrho=-\frac{p^2+p^2q-2q+4}{2(q-1)}\theta_3^2 -\frac 12 pq \theta_3 \circ \theta_4+(q-\frac 12 q^2)\theta_4^2.
\vspace{3pt}
\end{array}$ 
$\begin{array}{cl}
\vphantom{\displaystyle{A^{A^A}}} {\bf 1.3^1.6}: &  [u_1,u_3]= - u_2, \;   [u_1,u_4]=  u_1, \; [u_2,u_3]= u_1,  \; [u_2,u_4]= u_2, \;  [u_3,u_4]= e_1;\qquad  \varrho=2\theta_3^2 .
\vspace{3pt}
\end{array}$ 
\\$\begin{array}{cl}
\vphantom{\displaystyle{A^{A^A}}} {\bf 1.3^1.7}: & [u_1,u_3]=\frac{1}{1+p} e_1+ \frac{p}{1+p} u_1-\frac{1}{1+p} u_2, \; [u_1,u_4]=-\frac{1}{1+p} e_1+ \frac{1}{1+p} u_1+\frac{1}{1+p} u_2,\\[4pt] & [u_2,u_3]= -\frac{1}{1+p} e_1+ \frac{1}{1+p} u_1+\frac{1}{1+p} u_2, \; [u_2,u_4]=-\frac{p}{1+p} e_1+\frac{p}{1+p} u_1+\frac{1+2p}{1+p} u_2, \; p \neq -1; \\[4pt]
&  \varrho=-\frac{p-1}{2(p+1)}\theta_3^2 +\frac{p}{p+1}  \theta_3 \circ \theta_4-\frac{p}{p+1} \theta_4^2.
\vspace{3pt}
\end{array}$ 
\\$\begin{array}{cl}
\vphantom{\displaystyle{A^{A^A}}} {\bf 1.3^1.8}: &  [u_2,u_3]=  u_1,  \;  [u_2,u_4]= u_2, \;  [u_3,u_4]= -u_3;\qquad \varrho=-\frac{1}{2}\theta_4^2.
\vspace{3pt}
\end{array}$ 
\\$\begin{array}{cl}
\vphantom{\displaystyle{A^{A^A}}}{\bf 1.3^1.9}: &  [u_2,u_3]=  p u_1,  \; [u_2,u_4]=-p e_1+ (p + 1) u_2, \; [u_3,u_4]= - p u_3, \quad p\neq 0; \\[4pt]
& \varrho=(p -\frac 12 p^2 - \frac 12)\theta_4^2.
\vspace{3pt}
\\ {\bf 1.3^1.11}: &   [u_2,u_3]=  - u_1,  \; [u_2,u_4]= e_1, \;  [u_3,u_4]=e_1+u_3;
 \qquad \varrho=-2\theta_4^2.
\vspace{3pt}
\end{array}$ 
\\$\begin{array}{cl}
\vphantom{\displaystyle{A^{A^A}}} {\bf 1.3^1.12}: &  [u_1,u_4]= u_1,\; 
[u_2,u_3]=  q u_1,  \; [u_2,u_4]=-pq e_1 +(p+q) u_2, \; [u_3,u_4]= (1-q) u_3,  \\[4pt]
& \; p,q \in \R; \qquad \varrho=\frac{1-(p-q)^2}{2}\theta_4^2.
\vspace{3pt}
\end{array}$ 
\\$\begin{array}{cl}
\vphantom{\displaystyle{A^{A^A}}} {\bf 1.3^1.13}: &  [u_1,u_4]= u_1, \; \; [u_2,u_3]=  \frac 12 u_1, \;[u_2,u_4]=-\frac p2 e_1 +(p +\frac 12 )u_2, \; [u_3,u_4]=e_1+ \frac 12 u_3, \\[4pt]
& p\in \R; \qquad \varrho=\frac{3+4p-4p^2}{8}\theta_4^2.
\vspace{3pt}
\end{array}$
\\$\begin{array}{cl}
\vphantom{\displaystyle{A^{A^A}}} {\bf 1.3^1.14}: &  [u_1,u_4]= u_1, \; 
[u_2,u_3]=  (1-p) u_1,  \; [u_2,u_4]=p(p-1) e_1 +u_2, \; [u_3,u_4]= e_1 +p u_3,  \\[4pt]
& p \neq \frac 12; \qquad \varrho=2p(1-p)\theta_4^2.
\vspace{3pt}
\end{array}$ 
\\$\begin{array}{cl}
\vphantom{\displaystyle{A^{A^A}}} {\bf 1.3^1.19}: &  [u_1,u_4]= u_1,\;
[u_2,u_3]=  u_1, \; [u_2,u_4]=-e_1+ u_1+2 u_2; \qquad \varrho=\frac 12 \theta_4^2.
\vspace{3pt}
\end{array}$ 
\\$\begin{array}{cl}
\vphantom{\displaystyle{A^{A^A}}} {\bf 1.3^1.20}: &  [u_2,u_3]= u_1, \; [u_2,u_4]=- u_1+ u_2 , \; [u_3,u_4]=- u_3; \qquad \varrho=-\frac 12 \theta_4^2.
\vspace{3pt}
\end{array}$
\\$\begin{array}{cl}
\vphantom{\displaystyle{A^{A^A}}} {\bf 1.3^1.21}: &  [u_1,u_4]= u_1, \; 
[u_2,u_3]= p u_1,  \; [u_2,u_4]=-p e_1 +(1-p) u_1+ (1+p) u_2 ,  \\[4pt] & [u_3,u_4]=(1-p) u_3 \quad p\neq 1; \\[4pt] & \varrho=(p- \frac 12 p^2 ) \theta_4^2.
\vspace{3pt}
\end{array}$
\\$\begin{array}{cl}
\vphantom{\displaystyle{A^{A^A}}} {\bf 1.3^1.22}: &  [u_1,u_4]= u_1, \; 
[u_2,u_3]= \frac 12 u_1,  \; [u_2,u_4]=-\frac 12 e_1 +\frac 12 u_1+ \frac 32 u_2 , \; [u_3,u_4]=e_1 +\frac 12 u_3; \\[4pt] & \varrho=\frac{3}{8} \theta_4^2.
\vspace{3pt}
\end{array}$
\\$\begin{array}{cl}
\vphantom{\displaystyle{A^{A^A}}} {\bf 1.3^1.28}: &  [u_1,u_3]= 2u_1, \;  [u_1,u_4]= 2u_2, \;
[u_2,u_3]= u_2,  \; [u_2,u_4]=e_1-\frac 12 u_1, \; [u_3,u_4]=u_4; \\[4pt] & \varrho=-\frac 32 \theta_3^2-\frac 12 \theta_4^2.
\vspace{3pt}
\end{array}$ 
 }}

\medskip
{
{\small
${\bf 1.4^1)}$  $\quad [e_1,u_2]= u_1, \; [e_1,u_3]= u_2$,  so that the invariant metrics are given by
$$g=-a\theta_1\circ \theta_3+a\theta_2^2 +b\theta_3^2+c\theta_3 \circ \theta_4+d\theta_4^2, \quad ad \neq 0.$$
The possible cases are the following. 
\\$\begin{array}{cl}
\vphantom{\displaystyle{A^{A^A}}} {\bf 1.4^1.2}: & [u_1,u_4]=  p u_1, \; [u_2,u_4]=(p-1)  u_2, \; [u_3,u_4]= (p-2)u_3, \quad p \in \R; \\[4pt]
&  \varrho=\frac{3a(p-1)^2}{d} \theta_1\circ \theta_3-\frac{3a(p-1)^2}{d} \theta_2^2 -\frac{b(3p^2-9p+8)}{d}\theta_3^2 -\frac{3c(p-1)^2}{d} \theta_3 \circ \theta_4-3(p-1)^2\theta_4^2.
\vspace{3pt}
\end{array}$
\\$\begin{array}{cl}
\vphantom{\displaystyle{A^{A^A}}} {\bf 1.4^1.9}: &  [u_1,u_3]= u_1, \; [u_2,u_3]= r e_1+ u_2+u_4,  \; 
[u_3,u_4]= p u_4, \quad p,r \in \R; \\[4pt]
&  \varrho=-(r+\frac{d}{2a}+p^2+p)\theta_3^2.
\vspace{3pt}
\end{array}$
\\$\begin{array}{cl}
\vphantom{\displaystyle{A^{A^A}}} {\bf 1.4^1.10}: &  [u_1,u_3]= u_1, \; [u_2,u_3]= r e_1 +u_2, \; [u_3,u_4]= p u_4, \; p,r \in \R; \quad  \varrho=-(r+p^2+p)\theta_3^2.
\vspace{3pt}
\end{array}$
\\$\begin{array}{cl}
\vphantom{\displaystyle{A^{A^A}}} {\bf 1.4^1.11}: &  [u_1,u_3]= u_1, \; [u_2,u_3]= r e_1 + u_2+u_4, \; [u_3,u_4]= u_1- u_4, \; r \in \R; \quad \varrho=-(r+\frac{d}{2a})\theta_3^2.
\vspace{3pt}
\end{array}$
\\$\begin{array}{cl}
\vphantom{\displaystyle{A^{A^A}}} {\bf 1.4^1.12}: &  [u_1,u_3]= u_1, \; [u_2,u_3]= r e_1 + u_2,  \; 
[u_3,u_4]= u_1- u_4, \; r \in \R; \quad \varrho=-r\theta_3^2.
\vspace{3pt}
\end{array}$
\\$\begin{array}{l}
\vphantom{\displaystyle{A^{A^A}}} {\bf 1.4^1.15},{\bf 1.4^1.16},{\bf 1.4^1.17}: \quad [u_2,u_3]= \varepsilon e_1 + u_4,  \; [u_3,u_4]= u_1, \; \varepsilon=1,-1,0 ; \quad \varrho=-(\varepsilon+\frac{d}{2a})\theta_3^2.
\vspace{3pt}
\end{array}$
\\$\begin{array}{l}
\vphantom{\displaystyle{A^{A^A}}} {\bf 1.4^1.21},{\bf 1.4^1.22}: \quad [u_2,u_3]= \varepsilon e_1 ,  \; [u_3,u_4]= u_1, \; \varepsilon=1,-1 ; \quad \varrho=-\varepsilon\theta_3^2.
\vspace{3pt}
\end{array}$
}}

\medskip
{
{\small
${\bf 2.2^1)}$  $\quad [e_1,u_1]= u_1, \; [e_1,u_3]= -u_3, \; [e_2,u_2]= u_1, \; [e_2,u_3]= -u_4$,  so that the invariant metrics are given by
\vspace{-1pt}
$$g=a\theta_1\circ \theta_3+c\theta_2^2 +a \theta_2 \circ \theta_4+b\theta_3^2, \quad a \neq 0.$$

The possible cases are the following. 
\\$\begin{array}{cl}
\vphantom{\displaystyle{A^{A^A}}} {\bf 2.2^1.2}: & [u_1,u_2]= e_2,\; [u_1,u_3]= u_4, \; [u_2,u_3]= (p-1)u_3, \; [u_2,u_4]= p u_4, \; [u_1,u_4]=  p u_1, \\[4pt] &   \; [u_2,u_4]=(p-1)  u_2, \; [u_3,u_4]= (p-2)u_3, \quad p \in \R; \qquad  \varrho=\frac{p^2-4}{2} \theta_2^2.
\vspace{3pt}
\end{array}$
\\$\begin{array}{cl}
\vphantom{\displaystyle{A^{A^A}}} {\bf 2.2^1.3}: & [u_2,u_3]= u_3,  \; [u_2,u_4]= u_4;\qquad  \varrho=\frac{1}{2} \theta_2^2.
\vspace{3pt}
\end{array}$
}}

\medskip
{
{\small
${\bf 2.5^1)}$  $\quad [e_1,u_2]= u_1, \;  [e_1,u_3]= -u_4,\;  [e_2,u_3]= -u_2, \;   [e_2,u_4]= u_1$,  so that the invariant metrics are given by
$$g=a\theta_1\circ \theta_3 +a \theta_2 \circ \theta_4+b\theta_3^2, \quad a \neq 0.$$

The possible cases are the following. 
\\$\begin{array}{cl}
\vphantom{\displaystyle{A^{A^A}}} {\bf 2.5^1.3}: & [u_1,u_3]= u_1, \; [u_2,u_3]=e_1+r e_2 +(1-h)u_2,  \; [u_2,u_4]= hu_1, \\[4pt] & [u_3,u_4]= -(r+h)e_1+s e_2-(1+h)u_4, \quad h \geq 0 \; \text{if} \; s \neq 0; \\[4pt] &   \varrho=(2r+h-\frac{h^2}{2}) \theta_3^2.
\vspace{3pt}
\end{array}$
\\$\begin{array}{cl}
\vphantom{\displaystyle{A^{A^A}}} {\bf 2.5^1.4}: & [u_1,u_3]= u_1, \;  [u_2,u_3]= r e_2 +(1-h)u_2,  \;  [u_2,u_4]= hu_1, \\[4pt]&[u_3,u_4]= -(r+h)e_1-(1+h)u_4, \;  [u_2,u_3]= u_3,  \; [u_2,u_4]= u_4\quad h\geq 0, r\in \R;  \\[4pt] &  \varrho=(2r+h-\frac{h^2}{2})\theta_3^2.
\vspace{7pt}
\end{array}$
\\$\begin{array}{cl}
\vphantom{\displaystyle{A^{A^A}}} {\bf 2.5^2.2:} & [e_1,u_2]= u_1, \;   [e_1,u_3]= -u_2, \;  [e_2,u_3]= u_4, \;   [e_2,u_4]= -u_1, \\[4pt] &
[u_1,u_3]= u_1, \; [u_2,u_3]=(p+s)e_1+r e_2+ u_2-2r u_4,  \; [u_2,u_4]= 2r u_1, \\[4pt] &  [u_3,u_4]=-r e_1+(p-s)e_2 -2r u_2 -u_4,  \quad p\in \R, \; r,s \geq 0.\end{array}$\vspace{4pt}\\
The invariant metrics and their Ricci tensor are given by  
$$g= a\theta_1\circ \theta_3 +b\theta_2^2 +b\theta_3^2+a\theta_4^2, \; a\neq 0, \qquad \varrho= 2(p+r^2)\theta_3^3 .$$
\vspace{3pt}
\\$\begin{array}{cl}
\vphantom{\displaystyle{A^{A^A}}} {\bf 3.3^1.1:} & [e_1,u_2]= u_2, \,  [e_1,u_4]= -u_4, \, [e_2,u_2]= u_1, \,   [e_2,u_3]= -u_4, \, [e_3,u_3]= -u_2, \,   [e_3,u_4]= u_1 , \\[4pt] & 
[u_1,u_3]= u_1, \, [u_2,u_3]=p e_3+ u_2, \, [u_3,u_4]=-p e_2  -u_4, \quad p\in \R. \end{array}$\vspace{4pt}\\
The invariant metrics and their Ricci tensor are given by  
$$g= a\theta_1\circ \theta_3 +a\theta_2\circ\theta_4 +b\theta_3^2, \; a\neq 0, \qquad\varrho= 2p\theta_3^2 .$$
\vspace{3pt}
\\$\begin{array}{cl}
\vphantom{\displaystyle{A^{A^A}}} {\bf 3.3^2.1:} & [e_1,u_2]= u_4, \,  [e_1,u_4]= -u_2, \,  [e_2,u_2]= u_1, \, [e_2,u_3]= -u_2, \, [e_3,u_3]= u_4, \,  [e_3,u_4]= -u_1, \\[4pt] & 
[u_1,u_3]= u_1, \, [u_2,u_3]=p e_2+ u_2, \, [u_3,u_4]=p e_3 -u_4, \quad p\in \R. \end{array}$\vspace{4pt}\\
The invariant metrics and their Ricci tensor are given by 
$$g= a\theta_1\circ \theta_3 +a\theta_2^2 +b\theta_3^2+a\theta_4^2, \; a\neq 0, \qquad\varrho= 2p\theta_3^2 .$$}}\vspace{-6pt}
\vspace{-10pt}
{\begin{rem}  
{\em  Note that the expressions { of the isotropy in cases ${\bf 1.1^1}$,} ${\bf 1.1^2}$ and ${\bf 2.2^1}$ given in \cite{K} are more general { than the ones we reported above,} but only in the cases we listed they give rise to non-trivial homogeneous Ricci solitons. 
}
\end{rem}
}

 \vspace{-13pt}

\noindent
\section{Classification of four-dimensional homogeneous Ricci solitons}
\setcounter{equation}{0}

{Solving  equation $\eqref{solit}$ for every  four-dimensional homogeneous pseudo-Riemannian manifold  $M=G/H$ with non-trivial isotropy, described in \cite{K}, we can  show the following classification result.}

\begin{theorem}\label{4dinvsol}{{The only  four-dimensional  pseudo-Riemannian homogeneous spaces $M=G/H$ with non-trivial isotropy,  admitting a  non-trivial  homogeneous   pseudo-Riemannian Ricci soliton  $g$  satisfying   equation \eqref{solit}  for some vector field $X \in \m$ and $\lambda \in \R$, are listed in the following {\em Table~I}, where the checkmark \lq\lq \checkmark\rq\rq means that all homogeneous pseudo-Riemannian spaces of that given form are Ricci solitons.}}\end{theorem}

{\begin{center}
\begin{tabular}{|l|c||l|c||l|c|}
\hline
Case$\vphantom{A^{A^{A}}}$  & Condition & Case & Condition & Case & Condition \\
\hline
${\bf 1.1^1.1}$$\vphantom{A^{A^{A}}}$ & $b=0 \neq d$ & ${\bf 1.1^1.2}$  & $p\neq \frac 12$ & ${\bf 1.1^1.5}$ & \checkmark  \\
\hline
${\bf 1.1^2.1}$$\vphantom{A^{A^{A}}}$& \checkmark & ${\bf 1.1^2.2}$$\vphantom{A^{A^{A}}}$ & $p \neq 1$ & ${\bf 1.1^2.6,7}$  &\checkmark \\
\hline
${\bf 1.3^1.1}$$\vphantom{A^{A^{A}}}$ & \checkmark &
${\bf 1.3^1.2}$$\vphantom{A^{A^{A}}}$ & $p=0$ & ${\bf 1.3^1.3}$$\vphantom{A^{A^{A}}}$   & \checkmark \\
\hline ${\bf 1.3^1.5}$ & $q=0$ &
${\bf 1.3^1.6}$$\vphantom{A^{A^{A}}}$ & \checkmark  & ${\bf 1.3^1.7}$$\vphantom{A^{A^{A}}}$ & $p=0$ \\
\hline  ${\bf 1.3^1.8}$$\vphantom{A^{A^{A}}}$   & \checkmark  &
${\bf 1.3^1.9}$  &  $p \neq 1$ & ${\bf 1.3^1.11}$ & \checkmark \\
\hline ${\bf 1.3^1.12}$$\vphantom{A^{A^{A}}}$ & \checkmark & 
${\bf 1.3^1.13}$&$\begin{array}{l} p=0 \; \text{or} \\ 4p^2-4p-3 \neq 0\end{array}$  &  $\vphantom{A^{A^{A}}}$${\bf 1.3^1.14}$ & $p(p-1) \neq 0$ \\
\hline ${\bf 1.3^1.19}$$\vphantom{A^{A^{A}}}$ & \checkmark  & 
${\bf 1.3^1.20}$ & \checkmark &  ${\bf 1.3^1.21}$ & $p(p-2) \neq 0$ \\
\hline ${\bf 1.3^1.22}$$\vphantom{A^{A^{A}}}$ & \checkmark &
${\bf 1.3^1.28}$  &\checkmark  & ${\bf 1.4^1.2}$$\vphantom{A^{A^{A}}}$  & $p=1$ \\
\hline ${\bf 1.4^1.9}$& \checkmark &
${\bf 1.4^1.10}$$\vphantom{A^{A^{A}}}$ & $r+p(p+1) \neq 0$ &  ${\bf 1.4^1.12}$ & \checkmark \\
\hline
 ${\bf 1.4^1.15}$$\vphantom{A^{A^{A}}}$ & \checkmark & ${\bf 1.4^1.16}$  &\checkmark  & ${\bf 1.4^1.17}$ & \checkmark  \\
\hline
${\bf 1.4^1.21}$$\vphantom{A^{A^{A}}}$ & \checkmark & ${\bf 1.4^1.22}$ & \checkmark & ${\bf 2.2^1.2}$& $p \neq \pm 2$  \\
\hline
${\bf 2.2^1.3}$$\vphantom{A^{A^{A}}}$ & \checkmark & ${\bf 2.5^1.3}$  & $2r+h-\frac 12 h^2 \neq 0$ &${\bf 2.5^1.4}$ &  $2r+h-\frac 12 h^2 \neq 0$ \\
\hline
$\vphantom{A^{A^{A}}}$${\bf 2.5^2.2}$ & $p \neq -r^2$ & ${\bf 3.3^1.1}$$\vphantom{A^{A^{A}}}$ & $p\neq 0$ & ${\bf 3.3^2.1}$$\vphantom{A^{A^{A}}}$ & $p \neq 0$\\
\hline
\end{tabular}\nopagebreak \\ \nopagebreak {Table {I}: Homogeneous pseudo-Riemannian $4$-spaces admitting Ricci solitons} $\vphantom{\displaystyle\frac{a}{2}}$
\end{center}
}

\medskip\noindent
{\bf Proof.} The complete classification follows from a case-by-case argument. One first excludes the case where $[u_i,u_j]_{\mathfrak m}=0$ for all indices $i,j$. { In fact, in such a case, any vector field $X \in \m$ is Killing and so, equation \eqref{algsol} reduces to the Einstein equation $\varrho=\lambda g$.} In the remaining cases, we computed  the Levi-Civita connection and the Ricci tensor using \eqref{LC} and \eqref{ric} and solved \eqref{algsol}. When some solutions of \eqref{algsol} were found, we excluded the Einstein cases (trivial Ricci solitons) and 
checked the invariancy of the solution under the isotropy action. 

As an example, we report here the calculations for cases ${\bf 1.1^2.6}$ and ${\bf 1.1^2.7}$, where both invariant and not invariant solutions of \eqref{algsol} occur. So, let $M=G/H$ be a four-dimensional homogeneous space, such that the Lie algebra $\g$ and the isotropy subalgebra $\h$ are determined by conditions
$[e_1,u_1]= u_3$, $[e_1,u_3]= -u_1$, $[u_1,u_3]= \varepsilon e_1$ and  $[u_2,u_4]= u_2$. 

We first determine invariant pseudo-Riemannian metrics on $M$, which are the same for all the cases ${\bf 1.1^2}$ corresponding to non-trivial Ricci solitons. Starting from an arbitrary bilinear symmetric form $g$, the isotropy representation as described by $[e_1,u_1]= u_3$, $[e_1,u_3]= -u_1$ easily yields that $g$ is invariant (that is, ${\psi} (e_1)^t \circ g +g \circ {\psi}(e_1)=0$) if and only if $g$ is of the form given in \eqref{first}. We then apply \eqref{LC} to {compute} $\Lambda (u_i)$ for all indices $i=1,..,4$, so describing the Levi-Civita connection of $g$. We get
$$\Lambda (e_1)= \left(\begin{array}{cccc} 
0 & 0 & -1 & 0 \\  0 & 0 & 0 & 0 \\  1 & 0 & 0 & 0 \\ 0 & 0 & 0 & 0  
\end{array}\right)
$$
and
{
$$\begin{array}{ll} \Lambda (u_1)=0, & \Lambda (u_3)=0, \\[4pt]
\Lambda (u_2):=\left(\begin{array}{cccc}
0&0&0&0 \\  0& \frac{bd}{bc-d^2} & 0 & \frac{bc}{bc-d^2} \\  0&0&0&0 \\ 0&-\frac{b^2}{bc-d^2}& 0 & -\frac{bd}{bc-d^2}  
\end{array}\right), & 
\Lambda (u_4):=\left(\begin{array}{cccc}
0 & 0 &0& 0 \\  0& \frac{d^2}{bc-d^2} & 0 &\frac{cd}{bc-d^2} \\  0&0&0&0 \\ 0&-\frac{bd}{bc-d^2}& 0 & -\frac{d^2}{bc-d^2}  
\end{array}\right). \end{array}
$$}
The curvature and Ricci tensors can be now deduced from the above formulas by a direct calculation applying \eqref{Curv} and \eqref{ric}. In particular, the Ricci tensor has the form {
\begin{equation}\label{roproof}
\begin{array}{l}\varrho=\varepsilon\theta_1^2 -\frac{b^2}{bc-d^2}\theta_2^2-\frac{bd}{bc-d^2}\theta_2 \circ \theta_4+\varepsilon\theta_3^2 -\frac{bc}{bc-d^2}\theta_4^2.\end{array}
\end{equation}
Moreover, it is easily seen from \eqref{first} and \eqref{roproof} that the manifold $(M=G/H,g)$ is Einstein if and only if $\frac{\varepsilon}{a}=-\frac{b}{bc-d^2}$.

We now consider an arbitrary vector field $X= \sum _{k=1}^4  x_k u_k \in \m$ and a real constant $\lambda$ and write down \eqref{algsol}. We find that $X$ and $\lambda$ determine a Ricci soliton if and only if the components $x_k$ of $X$ with respect to $\{u_k \}$ and $\lambda$ satisfy

\begin{equation}\label{syssol}
\left\{\begin{array}{l}
\varepsilon-a \lambda =0, \\[4 pt]
b\left((bc-d^2)(\lambda-2x_4)+b\right)=0, \\[4 pt]
(bc-d^2)\left(b x_2-d x_4+d \lambda \right) +bd =0, \\[4 pt]
(bc-d^2)\left(2d x_2+c \lambda \right) +bc =0.
\end{array}\right.
\end{equation}
It must be noted that components $x_1,x_3$ do not appear in \eqref{syssol}. We {now  solve the system \eqref{syssol}}. Taking into account $a(bc-d^2) \neq 0$ and excluding the Eistein cases, we find $\lambda =\frac{\varepsilon}{a}$ (and so, the Ricci solitons are not steady), $b=0$ and
$$\begin{array}{l}  x_2 = -\frac{\varepsilon c}{2ad}, \quad  x_4=\frac{\varepsilon}{a} , \end{array}$$
for any real value of $x_1,x_3$. Thus, $X=x_1 u_1-\frac{\varepsilon c}{2ad}u_2+x_3 u_3 +\frac{\varepsilon}{a}u_4$. Finally, again by {$[e_1,u_1]= u_3, \; [e_1,u_3]= -u_1$} we see at once that $X \in \m$ is invariant if and only if $X \in {\rm Span}\{e_2,e_4\}$. Therefore, there exists a two-parameter family of non-trivial (steady) Ricci solitons on $(M=G/H,g)$ determined by vector fields $X \in \m$ of  the above form, but among them, only $X=-\frac{\varepsilon c}{2ad}u_2+\frac{\varepsilon}{a}u_4$ is invariant and so, determines a homogeneous Ricci soliton.

{{{By  similar computations we conclude that the following cases occurs: 

\medskip
 ${\bf 1.1^1.1:}$ $g$  with  $b=0\neq d$, $\lambda$ is arbitrary  and  
${ X=  \frac{4a^2+d^2- 2 c a^2 \lambda}{8 a^2 d} u_2 + \frac 12 \lambda u_4} \; {\text{(invariant)}}.$

\smallskip
${\bf 1.1^1.2:}$ the possible cases are the following:

\begin{enumerate}
\item $p =1$, $g$ with $b=0 \neq c$, $\lambda$ is arbitrary and $X=x_1 u_1-\frac{2d\lambda-1}{4c} u_2+\lambda u_4$  (invariant if $x_1=0$).

\vspace{3pt}
\item $p \neq \frac 12$, $g$ with $b=0\neq c$, $\lambda=0$ and $X=x_1 u_1+\frac{2p-1}{4pc} u_2$  (invariant if $x_1=0$).

\end{enumerate}
}

\smallskip
 ${\bf 1.1^1.5:}$ $g$ with $b=0$, $\lambda =-\frac 1a$ and 
$$\begin{array}{l}X=x_1 u_1+\frac{d}{2ac} u_2+x_3 u_3-\frac 1a u_4 \quad {\text{(invariant if} \; x_1=x_3=0).} \end{array}$$

\smallskip
{ ${\bf 1.1^2.1:}$ $g$  with $b=0$,  $\lambda$ is arbitrary  and }
$X=\frac{4a^2+d^2-2\lambda a^2 c}{8 a^2 d} u_2 +\frac 12 \lambda u_4 \; {\text{(invariant)}.}$

${\bf 1.1^2.2:}$
$p \neq 0,1$, $g$ with $b=0$, $\lambda=0$ and $X=\frac{p-1}{dp} u_2 \; \text{(invariant).}$

\smallskip

 ${\bf 1.1^2.6}$ and ${\bf 1.1^2.7}$: $g$ with  $b=0$, $\lambda =\frac{\varepsilon}{a}$ and
$$
\begin{array}{l}
X=x_1 u_1-\frac{\varepsilon c}{2ad}u_2+x_3 u_3 +\frac{\varepsilon}{a}u_4 \; {\text{(invariant if} \; x_1=x_3=0).}
\end{array}
$$

\vspace{-3mm}
${\bf 1.3^1.1}$: $g$ with $d =0 \neq c$, $\lambda = 0$  and  
$ X=-\frac{15c}{8a^2} u_2 \; {\text{(invariant)}.}$

${\bf 1.3^1.2:}$} { $p=0$, $g$  with}  $d = b \neq 0$, $\lambda =  - \frac{1}{2 b}$  { and
$
X = x_1 u_1 +\left(x_1 -\frac{b+c}{2ab}\right)  u_2  - \frac{1}{2b} u_3  - \frac{1}{2b} u_4.$

\smallskip
 ${\bf 1.3^1.3:}$ $g$ with}  $d = b \neq 0$, $ \lambda =  - \frac{1}{2 b}$  { and 
$
X = x_1 u_1 +\left(x_1 -\frac{b+c}{2ab}\right)  u_2  - \frac{1}{2b} u_3  - \frac{1}{2b} u_4.
$

\smallskip
{ ${\bf 1.3^1.5}$  with  one of the following conditions:}

\begin{enumerate}

\item  $p=q=0$, { $g$ with } $c=0$, $b \neq d$,  $\lambda =\frac{2}{b-d}$  and
$ X=-\frac{d\lambda}{2a}u_1+x_2 u_2+\lambda u_4.$

\smallskip
\item    $q=0 \neq p$,  $g$ is arbitrary, $\lambda =0$ and
$X=-\frac{p^2+4}{4a p}u_2 \; {\text{(invariant)}.}$ 

\end{enumerate}

\smallskip
${\bf 1.3^1.6:}$   $g$ satisfies $c=0\neq b-d$, $ \lambda=\frac{2}{b-d}$   and 
$X=-\frac{d}{a(b-d)}  u_1+x_2 u_2 +\frac{2}{b-d}  u_4.$

\smallskip
${\bf 1.3^1.7:}$  $p=0$, $g$ is arbitrary, $ \lambda=0$   and 
$X=-\frac{1}{4a}  u_2 \;{\text{(invariant)}.}$

\smallskip
 ${\bf 1.3^1.8:}$ $g$ with  $b =0$, $c \neq 0$,  $\lambda =0$   and
$X= x_1 u_1 + x_2  u_2 + \frac{1}{4 c}  u_3.$

\smallskip
${\bf 1.3^1.9:}$  $p \neq 1$, $g$ is arbitrary, $ \lambda=0$ and  
$X=x_1 u_1 - \frac{(p-1)^2b} {4ac}  u_2 +  \frac{(p- 1)^2}{4cp} u_3.$ 
 
\smallskip
${\bf 1.3^1.11:}$ $g$ { with $c \neq 0$,} $\lambda =0$ { and}
$X=x_1 u_1 - \frac{b} {2ac}  u_2 -\frac{1}{2c} u_3.$

\smallskip
{ ${\bf 1.3^1.12}$ with  one of the following conditions:}

\begin{enumerate} 
\item $(p-q)^2 \neq 1$, $g$ is arbitrary, $\lambda =0$  and 
$X=-\frac{(p-q)^2-1}{4a} u_1 \; {\text{(invariant)}.}$
\vspace{3pt}
\item    $p=q=0$,  $g$ is arbitrary, $\lambda =0$  and 
$X= -\frac {1}{4a}u_1 +x_2 u_2 \; {\text{(invariant)}.}$

\vspace{3pt}\item   $p=q=0$, { $g$ with $b=0$, $\lambda=a (\neq 0)$ and
$X= \frac {2ad-1}{4a}u_1 +x_2 u_2-a u_4 .$}

\vspace{3pt}\item  $p=0$, $q=\frac 12$, { $g$ with $c=0$, $\lambda=a (\neq 0)$ and
$ X=\frac{8ad-3}{16a}u_1 +x_2 u_2 -a u_3 .$
}
\end{enumerate}

\smallskip
{ ${\bf 1.3^1.13}$ with   one of the following conditions:}

\begin{enumerate}
\item { $4p^2-4p-3 \neq 0$,} $g$ is arbitrary, $\lambda =0$   and  
${ X= \frac{4p-4p^2+3}{16a}u_1 \; {\text{(invariant)}.}}$

\vspace{3pt}
\item   $p=0$, { $g$ with} $b=c=0$, { $\lambda$ is arbitrary}  and
$$\begin{array}{l}{ X= -\frac {8d\lambda +3}{16a}u_1 +x_2 u_2 +\lambda u_4  \quad {\text{(invariant if} \; \lambda=0).}}\end{array}$$ 

\end{enumerate}

\smallskip
{ ${\bf 1.3^1.14:}$}  $p(p-1) \neq 0${, $g$ is arbitrary}, $\lambda =0$  and 
${ X=-\frac{p(p-1)}{a} u_1 \; {\text{(invariant)}.}}$

\smallskip
{ ${\bf 1.3^1.19:}$ $g$ is arbitrary, $\lambda =0$ and} 
${ X=\frac{1}{4a} u_1 } \; {\text{(invariant)}.}$

\smallskip
{ ${\bf 1.3^1.20}$  with one of the following conditions:}
\begin{enumerate}

\item  { $g$ is arbitrary,} $ \lambda=0$ and 
$X= x_1 u_1 +\frac{1}{4a}u_2 \; {\text{(invariant)}.}$

\vspace{4pt}\item { $g$ with} $b=0$, $\lambda =0$ and
${  X=x_1 u_1 -\frac{4cx_3+1}{4a} u_2 +x_3 u_3} \; {\text{(invariant if} \; x_3=0).}$

\end{enumerate}

\smallskip
{ ${\bf 1.3^1.21:}$  { $p(p-2) \neq 0$,} $g$  is arbitrary, $\lambda =0$  and}
${ X=-\frac{p(p-2)}{4a} u_1 \; {\text{(invariant)}.}}$

\smallskip
{ ${\bf 1.3^1.22:}$ $g$ is arbitrary, $\lambda =0$ and} 
${ X=\frac{3}{16a} u_1 \; {\text{(invariant)}.}}$

\smallskip
 ${\bf 1.3^1.28:}$ $g$ with  $c=0$ and  $b-6d \neq 0$, $\lambda =\frac{1}{2(b-6d)}$  and
$X=x_1 u_1 +\frac{b-9d}{2a(b-6d)} u_2 - \frac{1}{2(b-6d)} u_3 .$

\medskip
{ ${\bf 1.4^1.2:}$ $p=1$, $g$ is arbitrary, $\lambda =0$  and} 
$X=-\frac 1d (\frac ca  u_1 +u_4) \; {\text{(invariant)}.}$

\smallskip
{ ${\bf 1.4^1.9}$ with  one of the following conditions:}

\begin{enumerate}
\item  $g$ with { $d+2ra+2ap(p+1) \neq 0$,} $\lambda =0$ and 
${ X= \frac{d+2ra+2ap(p+1)}{4a^2} u_1 \; {\text{(invariant)}.}}$

\vspace{4pt}\item  $p=0$, { $g$ is arbitrary,} $\lambda =0$  and 
${ X= \frac{2ra+d}{4a^2} u_1 +x_4 u_4 } \; {\text{(invariant)}.}$

\end{enumerate}

\smallskip
{ ${\bf 1.4^1.10:}$ { $r+p(p+1) \neq 0$,} $g$ is arbitrary,  $\lambda=0$  and}  
${  X=\frac{r+p(p+1)}{2a} u_1 } \; {\text{(invariant)}.}$

\smallskip
 ${\bf 1.4^1.12}:$ $r \neq 0$, $g$ is arbitrary, $\lambda =0$ and
$X=-\frac{r}{2a} u_1  \; {\text{(invariant)}}.$

\smallskip
{ ${\bf 1.4^1.15, 1.4^1.16}$ and ${\bf 1.4^1.17:}$ $2\varepsilon a +d \neq 0$, $g$ is arbitrary, $\lambda =0$ and}
$$\begin{array}{l}{ X=x_1 u_1 -\frac{2\varepsilon a +d}{4a^2} u_4 } \quad {\text{(invariant)}.}\end{array}$$

\smallskip
 ${\bf 1.4^1.21}$ and ${\bf 1.4^1.22:}$ $g$ is arbitrary, $\lambda =0$ and
$$\begin{array}{l}X=x_1 u_1 +x_2 u_2 -\frac{\varepsilon}{2a} u_4 \quad {\text{(invariant if} \; x_2=0).} \end{array}$$

\medskip
${\bf 2.2^1.2}:$   $p \neq \pm 2$, $g$ is arbitrary, $\lambda =0$ and 
$ X=- \frac{p^2-4}{4ap} u_4 \; {\text{(invariant)}.}$

\smallskip
 { ${\bf 2.2^1.3}$  with one of the following conditions:}

\begin{enumerate}

\item $g$  { is}  arbitrary, $\lambda =0$ and $X=x_1 u_1-\frac{1}{4a}u_4 \; {\text{(invariant if} \; x_1=0).}$

\item  { $g$ with} $b=0$,  { $\lambda$ is arbitrary}  and
$$\begin{array}{l}X=x_1 u_1-\lambda u_2 +\frac{2c\lambda -1}{4a}u_4 \quad {\text{(invariant if} \; \lambda=x_1=0).}\end{array}$$

\end{enumerate}

\medskip
 { ${\bf 2.5^1.3}$ and ${\bf 2.5^1.4:}$ { $2r+h-\frac 12 h^2 \neq 0$,} $g$ is arbitrary, $\lambda =0$ and }
$$\begin{array}{l}{ X=\frac{4r+2h-h^2}{4a} u_1} \quad {\text{(invariant)}.}\end{array}$$

\smallskip
 ${\bf 2.5^2.2:}$ $p \neq r$, $g$ is arbitrary, $\lambda =0$ and 
$X=\frac{p-r}{a} u_1 \; {\text{(invariant)}.}$

\smallskip
 ${\bf 3.3^1.1:}$ $p \neq 0$, $g$ is arbitrary, $\lambda =0$ and 
$X=\frac{p}{a} u_1  \; {\text{(invariant)}.}$

\smallskip
${\bf 3.3^2.1:}$ $p \neq 0$, $g$ is arbitrary, $\lambda =0$ and 
$X=\frac{p}{a} u_1 \; {\text{(invariant)}.}$}}}
 $\Box$

\medskip
It is easy to check that vector fields $X$ of any causal character occur in the list of vector fields determining homogeneous Ricci solitons {in the proof of} Theorem~\ref{4dinvsol}. In particular, several examples occur where $X$ is a light-like vector field.

\bigskip\noindent
\section{Geometric properties of four-dimensional homogeneous Ricci solitons}
\setcounter{equation}{0}

By Theorem~\ref{4dinvsol}, there exist a large number of four-dimensional pseudo-Riemannian homogeneous Ricci solitons with non-trivial isotropy. It is a natural problem to investigate the geometry of these examples, for instance whether they are gradient Ricci solitons, their curvature properties and their relationship with some other geometric structures. Readers interested in the geometry of the homogeneous Ricci soliton metrics  can easily work out several interesting cases from the above classification. A few remarkable behaviours are listed below.

\smallskip
 \subsection{Gradient Ricci solitons} We recall that a  gradient Ricci soliton is a pseudo-Riemannian manifold $(M,g)$, together with a smooth function $f$ on $M$, such that equation \eqref{solit} holds with {$X= Hess(f)$}. Several classification and rigidity results hold in the specific case of gradient Ricci solitons. Some examples may be found in \cite{BGG},\cite{MS},\cite{PW1},\cite{PW2}.  
 
 A gradient Ricci soliton $(M, g)$ is said to be {\em rigid} if it is isometric to a quotient of
$N  \times
R^k$, where
$N$ is an Einstein manifold and the potential function is a
generalization
of the
Gaussian soliton (i.e., $h =\frac{\lambda}{2} \| x \|^2$ on the
Euclidean factor).  Any homogeneous Riemannian gradient Ricci soliton is rigid  (\cite{PW1}) and any
homogeneous Lorentzian 
gradient Ricci soliton $(M, g)$  is rigid  if $\dim M \leq
4$ (\cite{BGGG}).  By \cite{CZ},   a non-reductive homogenous   gradient Ricci soliton  in dimension four is rigid,  but the neutral signature case for a general  homogeneous gradient Ricci soliton  is still open. 

If $X$ is (locally) a gradient, then the one-form $\omega_X$ dual to $X$ must be closed. Checking this condition for the Ricci solitons listed in Theorem~\ref{4dinvsol}, we obtain the following.

\begin{prop}\label{GRS}
Let $(M,g,X)$ be any of the homogeneous Ricci solitons listed in Theorem~{\em \ref{4dinvsol}}, where $(M=G/H,g)$ is a four-dimensional pseudo-Riemannian homogeneous space and $X \in \m$ a solution of \eqref{solit}. The one-form $\omega_X$ dual to $X$  is never closed in the following cases:

\smallskip\noindent
 ${\bf 1.1^2:}$ ${\bf 2}$;

\smallskip\noindent
 ${\bf 1.3^1:}$ ${\bf 1}$, ${\bf 6}$, ${\bf 7}$, ${\bf 8}$, ${\bf 11}$, ${\bf 14}$, ${\bf 19}$, ${\bf 20}$, ${\bf 21}$ and ${\bf 22}$;

\smallskip\noindent
 ${\bf 1.1^4:}$ ${\bf 10}$, ${\bf 12}$, ${\bf 15}$, ${\bf 16}$ and ${\bf 17}$;

\smallskip\noindent
 ${\bf 2.2^1.2};$ 

\smallskip\noindent
 ${\bf 2.5^1:}$ ${\bf 3}$ and ${\bf 4}$; 

\smallskip\noindent
 ${\bf 2.5^2.2};$

\smallskip\noindent
 ${\bf 3.3^1.1};$

\smallskip\noindent
 ${\bf 3.3^2.1}.$ 

\medskip\noindent
 In the remaining cases, $\omega_X$ is closed under some restrictions on either the invariant metric $g$ or the coefficients describing the Lie algebra.
\end{prop} 

\noindent
Thus, in all the cases listed in Proposition~\ref{GRS},  the homogeneous Ricci soliton cannot be a gradient one.

\smallskip
\subsection{Conformally flat examples}  In general a  pseudo-Riemannian manifold $(M,g)$ of dimension $n \geq 4$ is (locally) conformally flat if and only if
$$ W=R -\frac{1}{n-2} (\varrho -\frac{\tau}{n}g) \odot g -\frac{\tau}{2n(n-1)}g \odot g =0,$$
where $W$ is the Weyl curvature tensor, $\tau$ the scalar curvature of $(M,g)$ and $\odot$ denotes the Kulkarni-Nomizu product of two symmetric two-tensors. With regard to homogeneous Ricci soliton metrics listed in  Theorem~{\ref{4dinvsol}}, we have the following.

\begin{theorem}\label{cflat}
A four-dimensional homogeneous Ricci soliton $(M,g)$ (with non-trivial isotropy), as classified in Theorem~{\em \ref{4dinvsol}}, is conformally flat if and only if it corresponds to one of the following cases:

\smallskip
 ${\bf 1.1^1.1}$;

\smallskip
${\bf 1.1^2.1}$;  

\smallskip
${\bf 1.3^1.2}$;  ${\bf 1.3^1.5}, (1)$ with $2c=-pd$;  ${\bf 1.3^1.7}$ with $2c=-d$;  ${\bf 1.3^1.8}$;  ${\bf 1.3^1.9}$ with $p=-1$; ${\bf 1.3^1.19}$ with $b=0$; ${\bf 1.3^1.21}$ with either $p=\frac 12$ or $b=0$; ${\bf 1.3^1.28}$ with $b=2d$;

\smallskip
${\bf 1.4^1.2}$ with $b=0\neq d$; ${\bf 1.4^1.9}$ with $p=-1/2$ and $r=-\frac{a+4d}{4a}$;
${\bf 1.4^1.15}$ with $d=-a$;

\smallskip
${\bf 2.2^1.2}$; ${\bf 2.2^1.3}$;

\smallskip
${\bf 2.5^2.2}$ with $s=0$;

\smallskip
${\bf 3.3^1.1}$;

\smallskip
{${\bf 3.3^2.1}$.}
\end{theorem}

\noindent
A conformally flat Einstein pseudo-Riemannian manifold is of constant curvature and so, symmetric. Moreover, a conformally flat (locally) homogeneous {\em Riemannian} manifold is again (locally) symmetric \cite{Ta} (see \cite{C1} for some three-dimensional Lorentzian non-symmetric examples). 

However, replacing the Einstein condition with its generalization \eqref{solit}, these rigidity results do not hold any more. In fact, checking the (local) symmetry condition $\nabla R=0$ for the examples classified in Theorem~{\ref{cflat}}, we obtain the following.

\begin{cor}
Examples ${\bf 1.1^1.1}$, ${\bf 1.1^2.1}$,  ${\bf 1.3^1.5}$,  ${\bf 1.3^1.7}$, ${\bf 1.3^1.19}$, ${\bf 1.3^1.21}$,  ${\bf 1.3^1.28}$,  ${\bf 1.4^1.9}$,  ${\bf 2.2^1.2}$ with $b=0\neq p$, ${\bf 2.2^1.3}$, ${\bf 3.3^1.1}$ with $r \neq 0$ and ${\bf 3.3^2.1}$ with $r \neq 0$, are four-dimensional conformally flat pseudo-Riemannian homogeneous Ricci solitons which are not symmetric.
\end{cor}

\noindent
Lorentzian conformally flat gradient Ricci solitons were classified in \cite{BGG}. In particular, the ones determined by light-like vector fields are necessarily steady. Among four-dimensional pseudo-Riemannian homogeneous Ricci solitons with non-trivial isotropy, there are no examples which are at the same time Lorentzian (non-gradient), non-steady and defined by a light-like vector field $X \in \m$. However, there exist such examples with metrics of neutral signature $(2,2)$. For example, ${\bf 1.3^1.5}, (1)$ is non-steady and always defined by a light-like vector field, and is conformally flat when $c=-\frac 12 pd$.

\subsection{Self-dual and anti-self-dual examples} 

Let $(M,g)$ denote an oriented four-dimensional pseudo-Riemannian manifold of signature $(2,2)$ and $\Lambda^2$ the bundle of bi-vectors on $M$. The {\em Hodge star operator} $\ast:\Lambda^2\rightarrow\Lambda^2$ is the involution defined by
$$
\ast(e_1\wedge e_2)=e_3\wedge e_4,\ \ \ast(e_1\wedge e_3)=e_2\wedge e_4,\ \ \ast(e_1\wedge e_4)=-e_2\wedge e_3,
$$
where $\{e_1,e_2,e_3,e_4\}$ is a local oriented pseudo-orthonormal frame field of $M$, with $e_1, e_2$ space-like and $e_3 ,e_4$ time-like vector fields. Let $\Lambda^+$ and $\Lambda^-$ denote the subbundles of $\Lambda^2$ corresponding to the eigenvalues $\pm1$ of $\ast$, and $\mathcal R:\Lambda^2\rightarrow\Lambda^2$ the curvature operator defined from the curvature tensor $R$ of $g$ by
$$
g(\mathcal R(X\wedge Y),Z\wedge T)=g(R(X,Y)Z,T).
$$
Then, $\mathcal R$ admits an $SO_0 (2,2)$-irreducible decomposition
$$
\mathcal R=\frac{\tau}{12}Id_{\Lambda^2}+\mathcal B+
\left(
\begin{array}{cc}
\mathcal W^+ & 0 \\
0 & \mathcal W^- \\
\end{array}
\right),
$$
where $\tau$ is the scalar curvature, $\mathcal B$ the traceless Ricci tensor and $\W^\pm$ are the operators on $\Lambda^\pm$ induced by the Weyl conformal tensor $W$. Observe that the metric $g$ is Einstein exactly when $\mathcal B=0$, conformally flat when $\W=0$. The metric $g$ said to be {\em self-dual} (respectively, {\em anti-self-dual}) if $\W^+=0$ (respectively, $\W^-=0$).

A conformally flat metric is trivially both self- and anti-self-dual. Looking for (anti-)self-dual metrics, among the Ricci soliton metrics classified in Theorem~{\ref{4dinvsol}} but not included in Theorem~{\ref{cflat}}, we can find several examples. In particular:
\begin{itemize}
\item[(a)] ${\bf 1.3^1.1}$ {\em  is self-dual when $d(bd-c^2)=0$ and never anti-self-dual;}
\vspace{3pt}\item[(b)]${\bf 1.3^1.3}$ {\em is never self-dual and always anti-self-dual.}
\end{itemize}

\subsection{Existence of examples which are not solvsolitons} A solvsoliton is a Ricci soliton left-invariant metric $g$ on a solvable Lie group $G$. All known examples of homogeneous Riemannian Ricci soliton metrics on non-compact homogeneous manifolds are isometric to some solvsolitons (\cite{Jab}, Remark 1.5). On the other hand, some of the examples listed in Theorem~\ref{4dinvsol} show that {\em there exist homogeneous (and also invariant) pseudo-Riemannian Ricci solitons which are not isometric to solvmanifolds}. 
In fact, we found some homogeneous and invariant Ricci solitons in cases ${\bf 1.3^1.1}$ and ${\bf 1.4^1.2}$. It is easily seen that such homogeneous pseudo-Riemannian manifolds correspond respectively to cases ${\bf B1}$ and ${\bf A2}$ in the classification of {\em non-reductive} four-dimensional homogeneous pseudo-Riemannian manifolds obtained in \cite{Fels} (see also \cite{CF}). These spaces are diffeomorphic to $\mathbb R^4$ \cite{Fels} and so, non-compact. As they do not admit any reductive decomposition, obviously they cannot be isometric to any Lie group (solvable or not).

{\subsection{On the Segre type of the Ricci operator of the examples}
By definition, the {\em Ricci operator} $Q$ of a pseudo-Riemannian manifold $(M,g)$ is defined by $g(Q(X),Y)=\varrho(X,Y)$, for all tangent vector fields $X,Y$. Since the Ricci tensor $\varrho$ is symmetric, the Ricci operator $Q$ is self-adjoint. Hence,  if $g$ is Riemannian, then $Q$ is diagonalizable, as there always exists a local orthonormal basis of eigenvectors for $Q$. On the other hand, in general pseudo-Riemannian settings, a self-adjoint operator needs not be diagonalizable, and its different canonical forms ({\em Segre types}) are determined accordingly to the roots of the minimal polynomial.   For all the left-invariant Ricci solitons on three-dimensional Lorentzian Lie groups, the Ricci operator is not diagonalizable with a unique eigenvalue of multiplicity three \cite{isr}. A similar behaviour is found for Ricci soliton homogeneous metrics of non-reductive four-dimensional homogeneous manifolds \cite{CF}.

Following the standard terminology (see for example \cite[Section~5.1]{ste}, to which we may refer for more details), given a self-adjoint operator with respect to a nondegenerate inner product, its Segre type  lists between square brackets the sizes of Jordan blocks in the decomposition of the operator.

In particular, for the Ricci operator $Q$ of a four-dimensional Lorentzian manifold, the possible cases are the following:
\begin{enumerate}
	\item {\em Segre type $[111,1]$}: $Q$ is diagonalizable.
	\item {\em Segre type $[11, 1 \bar{1} ]$}: $Q$ has two real eigenvalues (which may coincide) and two complex conjugate eigenvalues.
	\item {\em Segre type $[11,2 ]$}: $Q$ has three real eigenvalues (some of which may coincide), one of which has multiplicity two and each associated to a one-dimensional eigenspace.
  \item {\em Segre type $[1,3]$}: the Ricci operator has two real eigenvalues (which may coincide), one of which has multiplicity three and each associated to a one-dimensional eigenspace.
\end{enumerate}
In the above list, the comma separates eigenvalues corresponding to space-like eigenvectors from those corresponding to time-like and light-like eigenvectors. In the different cases,  round brackets are used to group together different blocks referring to the same eigenvalue (when this occurs, in this context the self-adjoint operator is said to be \lq\lq degenerate\rq\rq, clearly not in the usual sense).

In the same way (see \cite{ste}), Segre types of the Ricci operator of a pseudo-Riemannian metric of neutral signature $(2,2)$ are classified into cases $[11,11]$, $[1,11\bar 1]$, $[1\bar 11\bar 1]$, $[1,12]$, $[22]$, $[21\bar 1]$, $[2\bar2]$, $[13]$, $[1,3]$ and $[4]$.
For instance, in the cases ${\bf 1.1^2.6,7}$, the homogeneous metric $g_{abcd}$ given by \eqref{first} is a Ricci soliton if and only if $b=0$. Equations \eqref{first} and \eqref{roproof} yield that the Ricci operator is then given by
$$Q=\left(\begin{array}{cccc} \frac{\varepsilon}{a} & 0 & 0 & 0 \\ 0 & -\frac{b}{bc-d^2} & 0 & 0 \\ 0 & 0 & \frac{\varepsilon}{a} & 0 \\ 0 &  0 & 0 & -\frac{b}{bc-d^2}\end{array}\right)$$
and thus, $Q$ is diagonalizable (more precisely, of Segre type [(11)(1,1)] when $g_{abcd}$ is Lorentzian, and [(11),(11)] when it is 
neutral). 

Clearly, the cases of Segre type $[(111,1)]$ and $[(11,11)]$ respectively  identify the invariant Einstein Lorentzian and neutral metrics, of four-dimensional pseudo-Riemannian manifolds, where $Q$ is diagonalizable and all Ricci eigenvalues coincide.  

In the cases under study, following a case-by-case argument, one can  determine by standard calculations the Segre types of the Ricci operator occurring for all examples of homogeneous Ricci solitons. However, this would be a long task to undertake in the present paper, and often the same family of invariant metrics contains examples with Ricci operators of different Segre types. For this reason, we limit ourselves to point out that Ricci solitons examples turn out to concern only a few non-diagonalizable Segre types of the Ricci operator: 
\begin{itemize}
\item[(a)] Segre type [(1,12)] (for example, in case ${\bf 2.2^1.3}$);
\vspace{3pt}\item[(b)] Segre type [(11,2)] (for example, in case ${\bf 1.3^1.8}$);
\vspace{3pt}\item[(c)] Segre type [(2,2)] (for example, in case ${\bf 1.3^1.28}$).
\end{itemize}}

\subsection{Examples which are not K\"ahler} In the Riemannian case, most of the known examples of non-trivial Ricci solitons correspond to K\"ahler metrics ({\em K\"ahler-Ricci solitons}), that is, they admit a K\"ahler structure and equation \eqref{solit} holds for a holomorphic vector field {X}  (see \cite{Cao, DW} and references therein). However, we have the following result.

\begin{theorem}
None of the examples listed in Theorem~{\em\ref{4dinvsol}} is (pseudo-)K\"ahler.  
\end{theorem}

\medskip\noindent
{\bf Proof.} Invariant pseudo-K\"ahler structures on four-dimensional homogeneous pseudo-Riemannian manifolds with non-trivial isotropy were classified by the present authors in \cite{CF2}. Comparing the classification of invariant K\"ahler metrics given in \cite{CF2} with the one of invariant Ricci soliton metrics listed in Theorem~\ref{4dinvsol}, it turns out that the only candidates as K\"ahler-Ricci solitons would be ${\bf 1.1^2:}$ ${\bf 1, 6,7}$.  However, in these cases, none of vector fields satisfying \eqref{solit} is holomorphic. 

As an example, we report the {computations}  for ${\bf 1.1^2.1}$. By Theorem~\ref{4dinvsol}, the invariant metric $g$ described in \eqref{first} is a Ricci soliton if and only if $b=0$. This is also a K\"ahler metric \cite{CF2}. In fact, $g$ is compatible with the two-parameter family of complex structures
$$
J_{\alpha \beta} = \pm \left(\begin{array}{cccc} 0 & 0 & 1 & 0 \\ 0 & \sqrt{-\alpha \beta -1} & 0 & \alpha \\ -1 & 0 & 0 & 0 \\ 0 & \beta & 0 & -\sqrt{-\alpha \beta -1}  \end{array}\right),
$$
for any real constants $\alpha, \beta$ such that $-\alpha \beta -1 \geq 0$, and the pair $(g, J_{\alpha \beta})$ is pseudo-K\"ahler if and only if  
$$
c (\alpha -\beta) + d \sqrt{-\alpha \beta -1} = 0.
$$
Again by Theorem~\ref{4dinvsol}, vector fields for which $g$ satisfies \eqref{solit} are given by 
$$X=\frac{4a^2+d^2-2\lambda a^2 c}{8 a^2 d} u_2 +\frac 12 \lambda u_4, $$ 
for any real value of $\lambda$. A standard calculation yields that $\mathcal L _X J_{\alpha \beta} \neq 0$. Therefore, $X$ is not holomorphic $\Box$

\smallskip
By the same argument used above, one can ask whether there exist homogeneous examples of {\em paraK\"ahler-Ricci solitons},  that is, invariant paraK\"ahler metrics, satisfying equation \eqref{solit} for a paraholomorphic vector field $X$. However, when we compare four-dimensional invariant paraK\"ahler metrics, as classified in \cite{CF2}, with the invariant Ricci soliton metrics listed in Theorem~\ref{4dinvsol}, we do not find invariant metrics which are both paraK\"ahler and also Ricci solitons. 

We remark that examples of invariant nontrivial Ricci solitons have also been found among four-dimensional paraK\"{a}hler Lie algebras 
(\cite{CHoust},\cite{Cproc}).

\end{document}